\documentclass[11pt]{amsart}
\usepackage{amsmath, amsthm, amssymb, latexsym,url}
\usepackage{hyperref}

\author{Joseph Vandehey}
\title{Absolutely abnormal, continued fraction normal numbers}
\date{\today}

\newtheorem{thm}{Theorem}[section]
\newtheorem{cor}[thm]{Corollary}

\allowdisplaybreaks

\begin{document}

\begin{abstract}
In this short note, we give a proof, conditional on the Generalized Riemann Hypothesis, that there exist numbers $x$ which are normal with respect to the continued fraction expansion but not to any base $b$ expansion. This partially answers a question of Bugeaud.
\end{abstract}

\maketitle

\section{Introduction}

A number $x\in [0,1)$ with base $b$ expansion $0.a_1a_2a_3\dots$ is said to be normal to base $b$ if for every finite string of digits $s=[d_1,d_2,\dots,d_k]$ with $d_i\in\{0,1,2,\dots,b-1\}$, we have
\[
\lim_{n\to\infty} \frac{A_{s,b}(n;x)}{n} = \frac{1}{b^k},
\]
where $A_{s,b}(n;x)$ denotes the number of times $s$ appears in the string $[a_1,a_2,\dots,a_n]$. Similarly, a number $x\in [0,1)$ with continued fraction expansion
\begin{equation}\label{eq:CFexpansion}
x= \cfrac{1}{a_1+\cfrac{1}{a_2+\cfrac{1}{a_3+\dots}}} = \langle a_1,a_2,a_3,\dots\rangle, \quad a_i \in \mathbb{N},
\end{equation}
is said to be normal with respect to the continued fraction expansion (or CF-normal) if for every finite string of positive integer digits $s=[d_1,d_2,\dots,d_k]$, we have
\[
\lim_{n\to \infty} \frac{B_s(n;x)}{n} = \mu(C_s),
\]
where $B_s(n;x)$ denotes the number of times the string $s$ appears in the string $[a_1,a_2,\dots,a_n]$, $\mu$ denotes the Gauss measure $\mu(A)= \int_A dx/(1+x)\log 2$, and $C_s$ denotes the cylinder set corresponding to $s$, the set of numbers in $[0,1)$ whose first $k$ continued fraction digits are the string $s$. 

If we let $\mathcal{N}_b\subset[0,1)$ denote the set of base-$b$ normal numbers and let $\mathcal{M}\subset[0,1)$ denote the set of CF-normal numbers, then by standard ergodic techniques all the sets $\mathcal{N}_b$, $b\ge 2$, and $\mathcal{M}$ have full Lebesgue measure. Despite the size of these sets, exhibiting any normal number is quite difficult. It is not known if $\pi$ is normal to any base $b$ or if it is CF-normal. All explicit examples of normal numbers (such as Champernowne's constant $0.12345678910111213\dots$ \cite{Champernowne}, which is normal to base $10$, or the Adler-Keane-Smorodinsky CF-normal number \cite{AKS}) were constructed to be normal. 

Our interest in this paper is understanding how the sets $\mathcal{N}_b$, $b\ge 2$, and $\mathcal{M}$ relate to one another. Steinhaus asked whether being $b$-normal implied normality for all other bases \cite{Steinhaus}. Maxfield showed that if $b,c\ge 2$ are positive integers such that $b^r=c$ for some rational number $r$, then $\mathcal{N}_b=\mathcal{N}_c$ \cite{Maxfield}. On the other hand, if $b^r=c$ and $r$ is not rational, then Cassels and Schmidt showed that $\mathcal{N}_b\setminus \mathcal{N}_c$ contains uncountably many points \cite{Cassels,Schmidt}. In fact, Schmidt showed that given any two sets of positive integers $S$ and $T$, the set of numbers in $\bigcup_{b\in S} \mathcal{N}_b \setminus \bigcup_{b\in T} \mathcal{N}_b$ is uncountably infinite, unless  prevented by Maxfield's theorem. 

Even fewer non-trivial examples of constructions of such selectively normal numbers exist. Bailey and Borwein showed that some of the Stoneham constants, which were constructed to be normal to a given base $b$, are not normal to certain other bases \cite{BB}. Greg Martin gave an explicit construction of an irrational number that is not in \emph{any} of the $\mathcal{N}_b$---such a number is said to be absolutely abnormal, as it is not normal to any base \cite{Martin}. 

Less is known about CF-normal numbers. Analogous definitions of normality exist for other types of continued fraction expansions besides the Regular Continued Fraction expansion used in \eqref{eq:CFexpansion}: for some of these expansions, it is known that the set of normal numbers for these expansions equals $\mathcal{M}$ (see, for instance Kraaikamp and Nakada \cite{KN} and the author \cite{VandeheyCF}).

In a similar vein to Steinhaus's question, Bugeaud has asked the following in the above terminology \cite[Problem~10.51]{Bugeaud}: Does there exist a $b\ge 2$ such that $\mathcal{M}\setminus \mathcal{N}_b$ is nonempty---that is, do there exist CF-normal numbers which are not normal to a given base $b$? Furthermore, is $\mathcal{M}\setminus \bigcup_{b=2}^\infty \mathcal{N}_b$ nonempty---that is, do there exist CF-normal but absolutely abnormal numbers? As with Cassels and Schmidt, we will show the answer is yes, but at present we can only offer a conditional yes.

\begin{thm}\label{thm:main}
On the Generalized Riemann Hypothesis, the set $\mathcal{M}\setminus \bigcup_{b=2}^\infty \mathcal{N}_b$ is dense and uncountable.
\end{thm}

We will use standard asymptotic notation in this paper. By $f(x)=O(g(x))$, we mean that there exists a constant $C$ such that $|f(x)|\le C |g(x)|$. By $f(x)=o(g(x))$, we mean that $\lim_{x\to\infty} f(x)/g(x) = 0$.

\section{Some additional facts}

\subsection{On continued fractions}

Let $x=\langle a_1,a_2,a_3,\dots\rangle$. We let $p_n/q_n$ denote the lowest terms expression for the truncated expansion $\langle a_1,a_2,\dots,a_n\rangle$ and refer to this as the $n$th convergent. When not clear which number $q_n$ refers to, we will state it as $q_n(x)$. The following facts can be found in most standard texts dealing with continued fractions, such as \cite{DK}.

The $q_n$'s satisfy a recurrence relation: with $q_{-1}=0$ and $q_0=1$ we have
\[
q_{n+1} = a_{n+1} q_n + q_{n-1}, \quad n=0,1,2,3,\cdots.
\]
We also  have that $q_n$ and $q_{n+1}$ are  relatively prime for $n\ge 0$. If $x$ is irrational, then its continued fraction expansion is infinite and unique. If $x$ is rational, then it is finite and not unique; however, the only finite expansions we will consider are convergents and so non-uniqueness will not be a concern in this paper.

Moreover, we can state in some way how well the convergents approximate $x$. In particular, we have
\begin{equation}\label{eq:approx}
\left| x- \frac{p_n}{q_n} \right| \le \frac{1}{a_{n+1} q_n^2}.
\end{equation}
Moreover, assuming $x$ is irrational, then $x-p_n/q_n$ is positive if and only if $n$ is even.

If we let $x=\langle a_1,a_2,\dots,a_n\rangle$ and $s=[a_1,a_2,\dots,a_n]$, then the cylinder set $C_s$ consists of the interval between $x $ and $\langle a_1,a_2,a_3,\dots,a_n+1\rangle$. (We will not care about whether endpoints are included.) Moreover, the length of the cylinder set decays to $0$ as $n$ tends to infinity.

\subsection{On Artin's Conjecture}

We reference the following result of Lenstra \cite[Theorem~8.3]{Lenstra}.

\begin{thm}\label{thm:Lenstra}
Assume the Generalized Riemann Hypothesis. Let $g,f,a$ be positive integers and let $P_{g,f,a}$ denote the set of primes $p$ such that $g$ is a primitive root modulo $p$ and $p\equiv a \pmod{f}$. Let $d$ be the discriminant of $\mathbb{Q}(\sqrt{g})$.

The set $P_{g,f,a}$ is finite if and only if one of the following conditions is true:
\begin{enumerate}
\item there exists a prime $q$ such that $q|f$, $a\equiv 1 \pmod{q}$, and $g$ is a perfect $q$th power;
\item we have $d|f$ and $\left( \frac{d}{a} \right) = 1$, where here we are using the standard Kronecker symbol; or,
\item we have $d|3f$, $3|d$,  
\[
\left( \frac{-d/3}{a} \right) = -1
\]
and $g$ is a perfect cube.
\end{enumerate}
\end{thm}

Conditional asymptotics on the size of $P_{g,f,a}$ can be found in the work of Pieter Moree \cite{PM1,PM2}.

From this we get the following simple corollary.

\begin{cor}\label{cor:Lenstra}
Assume the Generalized Riemann Hypothesis. Let $g,f,a$ be positive integers such that $f$ is relatively prime to both $g$ and $a-1$, and $g\ge 2$ is not a perfect square. Then $P_{g,f,a}$ is infinite.
\end{cor}

\begin{proof}
Let $g'$ be the largest squarefree divisor of $g$. First note that, by standard facts about discriminants, we have that 
\[
d=\begin{cases}
g', & \text{if } g'\equiv 1 \pmod{4},\\
4g', & \text{otherwise}.
\end{cases}
\]
 By our assumptions, $g'\ge 2$, so we must have that $d\ge 2$ and that $d$ shares a prime in common with $g$.

If $f$ is relatively prime to $a-1$, then for any prime $q|f$ we have that $a \neq 1 \pmod{q}$, so the first condition of Theorem \ref{thm:Lenstra} fails. If $f$ is relatively prime to $g$ then since $d$ shares a prime in common with $g$, we cannot have that $d|f$, and hence the second condition fails. By a similar argument, if $d|3f$ and $3|d$, then we must have that $d=3$; however, this is impossible since $d$ is either $0$ or $1$ modulo $4$.  Hence the third conditions fails as well, and so $P_{g,f,a}$ must be infinite.
\end{proof}

\section{Proof of the theorem}

Consider any $x\in[0,1)$ that is CF-normal. Let $x=\langle a_1,a_2,a_3,\dots\rangle$. Let $N$ be some large even integer, and let $\{X_i\}_{i=1}^\infty$ be the sequence of strings given by 
\[
X_1 = [a_1, a_2,\dots, a_N], \qquad X_i = [a_{2^{i-2}N+1}, a_{2^{i-2}N+2},\dots,a_{2^{i-1} N}], \text{ for }i\ge 2
\]
so that $x$ could be represented as the concatenation of the strings $X_i$. We will abuse notation slightly and allow ourselves to write $x=\langle X_1,X_2,X_3,\dots\rangle$.

We will consider a sequence of positive integers $L=\{\ell_1,\ell_2,\ell_3,\dots\}$ and create a new continued fraction expansion
\begin{equation}\label{eq:ydefn}
y=\langle X_1, \ell_1,\ell_2,\ell_3,\ell_4,X_2, \ell_5,\ell_6,\ell_7,\ell_8,X_3,\ell_9,\ell_{10},\ell_{11},\ell_{12},X_4,\dots\rangle= \langle a'_1,a'_2,a'_3,\dots\rangle.
\end{equation}
It is easy to see that $B_s(n;x)$ will differ from $B_s(n;y)$ by at most $O(\log_2 n)=o(n)$, and thus $y$ will be CF-normal regardless of what $L$ is. However, we claim that we can choose $L$ in such a way so that $y$ will be absolutely abnormal. Since $x$ and $N$ were arbitrary and $x$ and $y$ both belong to the cylinder set $C_s$ where $s=[a_1,a_2,\dots,a_N]$, this would also give us the density statement of the theorem. We will moreover show that for any fixed $x$ and $N$, there are uncountably infinite possibilities for $L$ such that $y$ is absolutely abnormal: If $x$ and $N$ are fixed, then each distinct choice of $L$ gives a distinct value of $y$, thus giving the uncountability statement of the theorem.

(If we instead wanted to prove a theorem about numbers that are absolutely abnormal and also \emph{not} CF-normal, then we could insert a sufficiently long string of $1$'s after each $X_i$. The rest of the construction would be unaltered.)

Let us illustrate how we intend to construct an $L$ that will leave $y$ absolutely abnormal. We will assume throughout the rest of this proof that $x$ and $N$ are fixed.

Let $\{b_1,b_2,b_3,\dots\}$ be a sequence of non-square bases $b_i\ge 2$ such that each base that is not a perfect square occurs infinitely often, for example $\{2,2,3,2,3,5,2,3,5,6,\dots\}$. Let  $n_i :=2^{i-1} N+4(i-1)$, so that $a'_{n_i}$ is the digit of $y$ that occurs at the end of the string $X_i$. This digit will then be followed by $a'_{n_i+j} = \ell_{4(i-1)+j}$ for $j=1,2,3,4$. Our choice of $\ell_{4(i-1)+1}$, $\ell_{4(i-1)+2}$, and $\ell_{4(i-1)+3}$ will then be made so that $q_{n_i+3}(y)$ is a power of $b_i$; that is, the $n_i+3$th convergent to $y$ will have a finite base $b_i$ expansion. Let $r_i\in \mathbb{Q}$ denote this $n_i+3$th convergent of $y$ and let $k_i$ denote the number of digits in its finite base $b_i$ expansion. 

Since we have assumed that $N$ is even, we have that $n_i+3$ is always odd and hence, for any choice of $L$, we have that $y-r_i$ is negative. Let us consider $r_i$ in its non-terminating base-$b_i$ expansion, i.e., the one that ends on the digit $b_i-1$ repeating. In particular, all digits after the $k_i$th digit will equal $b_i-1$.  By \eqref{eq:approx}, if $\ell_{4i}$ is sufficiently large---in particular, if $\ell_{4i}$ exceeds $b_i^{k_i^2}$---then the base-$b_i$ expansions of $y$ and $r_i$ only start to differ after the $k_i^2$th place. (Since $y<r_i$, this was why we considered the non-terminating expansion for $r_i$.) 

Suppose that we can find $L$ such that all of the $r_i$'s are rationals whose denominator is a power of $b_i$ and that $\ell_{4i}$ is always sufficiently large as defined in the previous paragraph. Then for any non-square base $b\ge 2$, there are infinitely many $b_i$'s such that $b_i=b$. For each such $i$, we know that at most $k_i$ of the first $k_i^2$ digits in the base-$b$ expansion of $y$ differ from $b-1$. Thus, if we let $s=[b-1]$, then 
\[ \limsup_{n\to \infty} \frac{A_{s,b}(n;y)}{n} =1.
\]
 So $y$ cannot be base-$b$ normal for any non-square base. For the remaining bases, we note that any square number $b\ge 2$ can be represented as a rational power of a non-square $c\ge 2$. By Maxfield's theorem as mentioned in the introduction, we know that $\mathcal{N}_c=\mathcal{N}_b$. Since we have already shown $y$ cannot be base-$b$ normal for any non-square base $b$, it must therefore not be base-$b$ normal for any square base $b$ either, and thus $y$ is absolutely abnormal.

It remains to show that we can find not just one, but uncountably many $L$'s with the desired properties. We will do this by showing that regardless of what $\{\ell_1,\ell_2, \dots, \ell_{4(i-1)}\}$ are, we can choose $\ell_{4(i-1)+1}, \ell_{4(i-1)+2}, \ell_{4(i-1)+3}$ in such a way so that $r_i$ has a denominator that is a power of $b_i$. Thus we can construct uncountably many $L$'s with the desired properties by first finding a desired $\ell_1, \ell_2,\ell_3$, then noting that we have countably many possibilities for $\ell_4$, namely all the values that are sufficiently large; and for each choice of $\ell_4$, we can find a desired $\ell_5, \ell_6, \ell_7$, and then have countably many possibilities for $\ell_8$; and so on.

Fix $i$ for now, and suppose that $\ell_1, \ell_2, \dots, \ell_{4(i-1)}$ have already been chosen. We will let $q_n$ refer to the denominator of the rational number $\langle a'_1,a'_2,\dots,a'_n\rangle$, where the $a'_j$'s are as in \eqref{eq:ydefn}, presuming we have chosen the values of enough $\ell_j$'s for this definition to make sense.

Since $q_n$ and $q_{n+1}$ are always relatively prime for any $n$, we know that for any prime $p$, there exists an $\ell\in \mathbb{N}$ such that
\[
\ell q_{n_i}+q_{n_i-1} \not\equiv 0 \pmod{p}.
\]
By the Chinese Remainder Theorem, we can thus find an $\ell \in \mathbb{N}$ such that $\ell q_{n_i} + q_{n_i-1}$ is relatively prime to both $b_i$ and $q_{n_i}-1$. We let $\ell_{4(i-1)+1}=\ell$ so that $q_{n_i+1}$ is now relatively prime to $b_i$ and $q_{n_i}-1$.

By our assumption that the Generalized Riemann Hypothesis is true, Corollary \ref{cor:Lenstra}, and our choice of $q_{n_i+1}$, there exist infinitely many primes congruent to $q_{n_i}$ modulo $q_{n_i+1}$ for which $b_i$ is a primitive root. We choose $\ell_{4(i-1)+2}\in \mathbb{N}$ so that $q_{n_i+2}=\ell_{4(i-1)+2}q_{n_i+1}+q_{n_i}$ is one of these primes.

Now since $q_{n_i+2}$ is a prime for which $b_i$ is a primitive root and since $q_{n_i+1}$ and $q_{n_i+2}$ are relatively prime, the following equation has a solution:
\[
q_{n_i+1} \equiv b_i^k \pmod{q_{n_i+2}}.
\]
We may suppose, without loss of generality that $k$ is a sufficiently large integer so that $b_i^k >2q_{n_i+2}$. Then we choose $\ell_{4(i-1)+3}$ so that $b_i^k = q_{n_i+3} =\ell_{4(i-1)+3} q_{n_i+2}+q_{n_i+1}$. But this shows that $q_{n_i+3}$ is a power of $b_i$ as desired and the proof is complete.

\section{Further questions}

It would be nice if we could give an unconditional proof of Theorem \ref{thm:main}. It would be especially nice if such a proof came as a result of a proof of the Generalized Riemann Hypothesis, but given the difficulty of such a problem, we propose the following question, which could lead to an unconditional result. Let $b\ge 2$ be fixed. For any rational number $p/q=\langle a_1,a_2,\dots,a_n\rangle$, does there exists a finite string $s=[a'_1,a'_2,\dots,a'_m]$ of length $m=o(n)$ such that the rational number $p'/q' = \langle a_1,a_2,\dots,a_n ,a'_1,a'_2,\dots,a'_m\rangle $ has a denominator that is a power of $b$? Can one replace $o(n)$ with $O(1)$?

We showed, assuming the Generalized Riemann Hypothesis, that $\mathcal{M}\setminus \bigcup_{b=2}^\infty \mathcal{N}_b$  is uncountable and dense. Bill Mance communicated the following question: what is the Hausdorff dimension of this set? Mance has studied the dimensional analysis of such differences sets in other cases, most notably for $Q$-Cantor series, and found that these difference sets often have full Hausdorff dimension \cite{ManceNT}.

\section{Acknowledgments}

The author acknowledges assistance from the Research and Training Group grant DMS-1344994 funded by the National Science Foundation. The author would also like to thank Kevin Ford, Greg Martin, and Paul Pollack for their helpful comments.

\end{document}